\documentclass[11pt]{article}

\usepackage{amsmath}
\usepackage{amssymb}
\usepackage{amsthm}

\def\fpd#1#2{{\displaystyle\frac{\partial #1}{\partial #2}}}
\def\spd#1#2#3{{\displaystyle\frac{\partial^2 #1}
{\partial #2\partial #3}}}

\def\vf#1{{\displaystyle\frac{\partial }{\partial #1}}}

\def\onehalf{{\textstyle\frac12}}

\def\supp{\mathop\mathrm{supp}}

\def\conn#1#2#3{\setbox1=\hbox{$\scriptstyle{#2}{#3}$}%
\setbox2=\hbox to\wd1{$\hfil\scriptstyle{#1}\hfil$}
\Gamma^{\!\box2}_{\!\box1}}

\def\TMO{T^\circ\!M}
\def\TxMO{T_x^\circ\!M}

\newcommand{\R}{\mathbb{R}}

\newcommand{\To}{T^\circ}

\newcommand{\TVO}{\To\kern-0.08em(\VM)}

\newcommand{\V}{\mathcal{V}}
\newcommand{\VM}{\V\kern-0.08em M}

\newtheorem{thm}{Theorem}
\newtheorem{lem}{Lemma}
\newtheorem{prop}{Proposition}
\newtheorem{cor}{Corollary}

\begin{document}
\title{The multiplier approach to the projective Finsler metrizability problem}

\author{M.\ Crampin, T.\ Mestdag\\
Department of Mathematics, Ghent University\\
Krijgslaan 281, B--9000 Gent, Belgium\\
and D.\,J.\ Saunders\\
Department of Mathematics, The University of Ostrava\\
30.\ dubna 22, 701 03 Ostrava, Czech Republic}

\maketitle

\begin{abstract}\noindent
This paper is concerned with the problem of determining whether
a projective-equivalence class of sprays is the geodesic class
of a Finsler function. We address both the local and the global
aspects of this problem. We present our results entirely in
terms of a multiplier, that is, a type $(0,2)$ tensor field
along the tangent bundle projection. In the course of the
analysis we consider several related issues of interest
including the positivity and strong convexity of
positively-homogeneous functions, the relation to the so-called
Rapcs\'ak conditions, some peculiarities of the two-dimensional
case, and geodesic convexity for sprays.
\end{abstract}

\subsubsection*{MSC}
53C60

\subsubsection*{Keywords}
Spray; projective equivalence; Finsler function; projective metrizability;
strong convexity; multiplier; geodesic convexity.

\section{Introduction}

A Finsler function --- a smooth function on the slit tangent
bundle $\tau:\TMO\to M$ of a smooth paracompact manifold $M$,
which is positive, positively (but not necessarily absolutely)
homogeneous, and strongly convex --- determines through its
Euler-Lagrange equations an equivalence class of sprays whose
base integral curves are its geodesics, where the curves with
given initial point and direction defined by two different
members of the class are the same up to an
orientation-preserving reparametrization. The geodesics of the
class, in other words, define oriented paths in $M$. Such sprays
are projectively equivalent. That is to say, if $\Gamma$ is one
such spray --- a vector field on $\TMO$ such that
$\tau_*\Gamma_{(z,y)}=y$ for any $x\in M$ and $y\in T_xM$,
$y\neq 0$, and such that $[\Delta,\Gamma]=\Gamma$ where $\Delta$
is the Liouville field --- then any other takes the form
$\Gamma-2P\Delta$ for some function $P$ on $\TMO$ which
satisfies $\Delta(P)=P$ (the choice of numerical coefficient is
made for later convenience). A set of sprays related in this way
is called a projective-equivalence class, or simply projective
class, of sprays.

The projective metrizability problem is the corresponding
inverse problem:\ given a projective class of sprays, does it
come from a Finsler function in the way just described? There
have been several publications on this question in the last few
years:\ see for example \cite{AP,Bacso,BM,Japan,Szi,SzVa}. It is
a curious fact that none of these uses what is, in the context
of the inverse problem of the calculus of variations in general,
probably the most studied, and certainly the historically most
significant, approach, namely that of the multiplier. The
necessary and sufficient conditions for the existence of a
Lagrangian for a given system of second-order ordinary
differential equations, when expressed as conditions on a
multiplier matrix, are known as the Helmholtz, or sometimes
Helmholtz-Sonin, conditions. Such conditions for systems with
two degrees of freedom were originally formulated by Douglas in
1941 in \cite{Doug}. The subject was resurrected by Henneaux,
and a more modern version of the Helmhotz conditions was given
by Sarlet, both in 1982 (see \cite{Hen} and \cite{Sar}). For a
recent survey see \cite{KP}.

Our first aim in this paper is to formulate the necessary and
sufficient conditions on a multiplier, conceived as a symmetric
type $(0,2)$ tensor field along the tangent bundle projection
$\tau$, for it to be the Hessian (with respect to the fibre
coordinates) of a positively-homogeneous function on $\TMO$ for
which a given projective class of sprays satisfies
the Euler-Lagrange equations. These necessary and sufficient
conditions play the role of Helmholtz conditions for the
projective Finsler metrizability problem. They differ from the
Helmholtz conditions for the general inverse problem of the
calculus of variations in several respects, arising from the
homogeneity of the sought-for Finsler function and the fact that
we work with a projective class of sprays rather than a single differential
equation field.

We must make a point of clarification here. As well as the
projective metrizability problem, there is another, related but
different and more straightforward, inverse problem in Finsler
geometry, that of determining whether a single given spray is
the canonical geodesic spray of a yet unknown Finsler function.
Given a Finsler function $F$ and its corresponding energy
$E=\onehalf F^2$, there are actually two candidates for the role
of a multiplier, namely the Hessians $h$ and $g$ of $F$ and $E$
respectively; they have significantly different properties in terms of
homogeneity and regularity. The Helmholtz conditions for the
canonical spray problem have been formulated previously (see
\cite{KrupkaS, Sar2}); the multiplier in this case is of the
type of $g$. By constrast we use the Finsler function rather
than the energy, and a multiplier of the type of $h$ rather than
$g$. We work with a projective class of sprays, and
avoid reference to the canonical spray; moreover we express our
results as far as possible in projectively-invariant terms.

We should also mention that when referring to a multiplier, and
indeed more generally, we usually abbreviate the expression
`tensor field along the tangent bundle projection' to simply `tensor field' or
`tensor'.

In addition to formulating the Helmholtz-like conditions for the
projective metrizability problem, we analyse in detail the
requirements on a Finsler function that it takes only positive
values and is strictly convex.

Our first results on the metrizability problem, though
global with respect to the fibres of $\tau$, or y-global, are
local with respect to $M$. As a third topic, we address the
question of how these results may be made global in $M$. In the
literature one will find mostly local results, but global
questions are barely touched upon. We find that there is a
cohomological obstruction on $M$ to combining locally-defined
Finsler functions into a global function; when this obstruction
vanishes the resulting function satisfies all of the
requirements of a Finsler function except positivity, and can be
turned into a genuine Finsler function (by the addition of a
total derivative) only locally. A well-known example of a spray
which exhibits precisely this behaviour is Shen's circle example
from \cite{Shen}.

Some authors, such as \'{A}lvarez Paiva (in \cite{AP}), deal
only with reversible paths, that is, paths which have no
preferred orientation; we, on the other hand, cover the more
general case of oriented paths, or sprays in the fully general
sense. At the level of Finsler geometry this distinction
corresponds to that between positively-homogeneous Finsler
functions (the general case) and absolutely-homogeneous
functions (the case discussed by \'{A}lvarez Paiva). We do
however specialise to the reversible case where appropriate, and
we begin the discussion of this in Section 2.

One of the difficulties of working directly with the Finsler
function, rather than with the energy, is that the conditions of
positivity and strong convexity required of a Finsler function
become somewhat tricky to deal with. We discuss this issue in
Section 3. The argument is based on the theorem from \cite{Bao}
in which the triangle inequality and the fundamental inequality
are established for Finsler functions:\ we prove a version which
doesn't assume ab initio that the function in question is
non-negative, and we use our result to show that any
positively-homogeneous function whose Hessian is positive
quasi-definite --- as nearly positive definite as it can be ---
can be turned into a Finsler function locally by the addition of
a total derivative.

Sections 2 and 3 deal with preliminary matters. In Section 4 we
turn to the projective metrizability problem, the main business
of the paper, and discuss the Rapcs\'{a}k conditions, which are
essentially ways of formulating or reformulating the
Euler-Lagrange equations. In Sections 5 and 6 we give the
necessary and sufficient conditions for the existence of a
Finsler function in terms of a multiplier. We discuss them
initially in local terms:\ local, that is, relative to $M$. In
fact throughout the paper, with the exception of Section 7, we
deal only with objects which are defined and smooth as functions
of $y^i$ for all $y\neq 0$, where $y=(y^i)$ are the fibre
coordinates on $TM$. That is to say, we deal almost entirely
with objects which are y-global. We discuss in Section 6 the
problem of extending our local results to results which are
global in $M$.

In the main body of the paper we make the assumption that $\dim
M\geq 3$. The two-dimensional case has some special features
which mean that our general results do not always apply. We make
some remarks about the two-dimensional case in Section 7. The
paper proper ends with the discussion of an illustrative
example, and some concluding remarks.

Most of our local calculations are carried out in coordinates.
We denote coordinates on $M$ by $(x^i)$, and the corresponding
canonical coordinates on $TM$ by $(x^i,y^i)$. We denote a spray
$\Gamma$ by
\[
\Gamma=y^i\vf{x^i}-2\Gamma^i\vf{y^i}.
\]
The corresponding horizontal local vector fields are denoted by
\[
H_i=\vf{x^i}-\Gamma^j_i\vf{y^j},\quad
\Gamma^i_j=\fpd{\Gamma^i}{y^j}.
\]
We write $V_i$ for $\partial/\partial y^i$
when it is convenient to do so.

The paper is written under the assumption that the reader is
familiar with the differential geometry of sprays. However, we
have given a number of basic results in an appendix, for
reference and to fix notations. As a general reference for
Finsler geometry we use the well-known book by Bao, Chern and
Shen, \cite{Bao}. The standard modern reference for the geometry
of sprays is Shen's book, \cite{Shen}; for a survey using a
coordinate-independent formalism see the paper by
B\'{a}cs\'{o} and Z.\ Szilasi, \cite{Bacso}.

A second appendix contains a discussion of the existence of
geodesically convex sets for a spray, which can be used in the
construction of so-called good open coverings, needed for the
global results in Section~6. While this material on convex sets,
which is of some interest in its own right, is not strictly
speaking new, it seems to have been lost from view for some
considerable time.

\section{Reversible sprays}

The geodesics of a spray are not in general reversible,
even as paths. In this brief section we consider the special class of
sprays for which the geodesics are reversible.

We denote by $\rho$ the reflection map of $\TMO$, namely
$(x,y)\mapsto(x,-y)$. For any spray $\Gamma$ the vector field
$\bar{\Gamma}=-\rho_*\Gamma$ is also a spray, which we call the
reverse of $\Gamma$. We say that $\Gamma$ is {\em reversible\/}
if it is projectively equivalent to its reverse, and {\em
strictly reversible\/} if the two are equal. Then the geodesic
paths of $\Gamma$ are reversible if and only if $\Gamma$ is
reversible. Moreover, the geodesics are reversible
as parametrized curves, which is to say that if $\gamma$ is a
geodesic so is $t\mapsto\gamma(-t)$, if and only if $\Gamma$ is
strictly reversible.

If a spray is reversible so are all sprays projectively
equivalent to it:\ that is to say, reversibility is a property
of a projective class. We shall show that the projective class
of a reversible spray contains a strictly-reversible spray.

\begin{prop}\label{rev}
If $\Gamma$ is reversible then there is a projectively-equivalent
spray which is strictly reversible.
\end{prop}

\begin{proof}

Let $\bar{\Gamma}$ be the reverse of $\Gamma$, and set
\[
\tilde{\Gamma}=\onehalf(\Gamma+\bar{\Gamma}).
\]
Then $\tilde{\Gamma}$ is a spray, and
\[
\rho_*\tilde{\Gamma}=\onehalf(\rho_*\Gamma+\rho_*\bar{\Gamma})
=-\onehalf(\bar{\Gamma}+\Gamma)=-\tilde{\Gamma}.
\]
Thus $\tilde{\Gamma}$ is strictly reversible. But $\bar{\Gamma}$
is projectively equivalent to $\Gamma$, from which it follows
immediately that $\tilde{\Gamma}$ is projectively equivalent to
$\Gamma$ also.
\end{proof}

\section{Positivity and strong convexity}

For a positively-homogeneous smooth function $F$ on $\TMO$ to be
a Finsler function it must be positive and strongly convex:\
that is to say,
\[
g_{ij}=F\spd{F}{y^i}{y^j}+\fpd{F}{y^i}\fpd{F}{y^j}
\]
must be positive definite. The conditions for a
projective class of sprays to be metrizable to be discussed
below consist in the first instance of conditions for the
existence of a positively-homogeneous  function $F$ whose
Euler-Lagrange equations are satisfied by the given sprays,
which we may describe as differential conditions. These
differential conditions must then be supplemented by algebraic
conditions which will ensure that there is such a function $F$
which is positive and strongly convex. Formulating such
algebraic conditions is somewhat tricky, for reasons we now
explain.

It follows from the fact that $F$ is homogeneous of degree 1 in
the fibre variables that its Hessian satisfies
\[
y^j\spd{F}{y^i}{y^j}(y)=0.
\]
At any point $(x,y)\in\TMO$ such that $F(x,y)\neq 0$, any vector
$u$ may be written uniquely as the sum of a scalar multiple of
$y$ and a vector $v$ which is $g$-orthogonal to $y$, namely
\[
v=u-\lambda y\quad\mbox{where }
\lambda=\frac{1}{F}u^k\fpd{F}{y^k}.
\]
Then
\[
g_{ij}u^iu^j=g_{ij}v^iv^j+\lambda^2g_{ij}y^iy^j=
F\spd{F}{y^i}{y^j}v^iv^j+\left(u^k\fpd{F}{y^k}\right)^2.
\]
Thus provided that $F(x,y)$ is positive, $g$ is positive
definite if and only if the Hessian is {\em positive
quasi-definite}, in the sense that
\[
\spd{F}{y^i}{y^j}(y)v^iv^j\geq0
\]
with equality only if $v$ is a scalar multiple of $y$.

So for $F$ to be a Finsler function it is necessary that its
Hessian be positive quasi-definite; but this property is not
sufficient. The problem is that we require $F$ to be positive,
but there is no way of ensuring this by imposing a condition on
the Hessian. Note that if $F$ is positively homogeneous so is
any function obtained by adding a total derivative to it, that
is, any function of the form
\[
F+y^i\fpd{\phi}{x^i}
\]
where $\phi$ is a function on $M$; the new function has the same
Hessian, and its Euler-Lagrange equations are satisfied by the
given sprays if those of $F$ are. But adding a total derivative
to a Finsler function may destroy the Finsler property.

Conversely, however, given a positively-homogeneous function $F$
whose Hessian is positive quasi-definite in the sense specified
above, it is always possible to modify $F$ by the addition of a
total derivative so as to obtain a local Finsler function, that
is, a function on $T^\circ U$ for some open neighbourhood $U$ of
any point in $M$, which is positive, positively homogeneous and
strongly convex. We shall prove this below.

It turns out, as will become apparent below, that there is no
difficulty about positivity if $F$ is absolutely rather than
just positively homogeneous, that is if $F(ky)=|k|F(y)$ for all
$k\in\R$.

The results in this section are not to be confused with that of
Lovas \cite{Lov}, who showed in effect that assuming positivity
one can prove strong convexity:\ he showed, that is, that for a
positive, positively-homogeneous function $F$, if
\[
\spd{F}{y^i}{y^j}(y)v^j=0
\]
only when $v$ is a scalar multiple of $y$ then the Hessian is
positive quasi-definite (this formulation is to be
found in \cite{hil4}). Lovas's result is converse to ours.

For the first part of the argument we fix $x\in M$, that is, we
essentially work in $\R^n$. The following result is Theorem
1.2.2 of \cite{Bao}, but without the assumption that $F$ is
non-negative.
\begin{prop}\label{tri}
Let $F$ be $C^\infty$ on $\R^n-\{0\}$, continuous on $\R^n$,
positively homogeneous, and suppose that its Hessian is
non-negative in the sense that for any $w\in\R^n$
\[
w^iw^j\spd{F}{y^i}{y^j}(y)\geq 0.
\]
Then for all $y_1,y_2\in\R^n$
\[
F(y_1)+F(y_2)\geq F(y_1+y_2)\quad\mbox{(the triangle inequality);}
\]
and for all $y,z\in\R^n$ with $y\neq 0$
\[
 F(z)\geq z^i\fpd{F}{y^i}(y)\quad\mbox{(the fundamental inequality).}
\]
\end{prop}

\begin{proof}
Provided that $0$ does not lie in the line segment $[y,y+w]$ we
can apply the second mean value theorem to the function
$t\mapsto F(y+tw)$ to obtain
\[
F(y+w)=F(y)+w^i\fpd{F}{y^i}(y)
+w^iw^j\spd{F}{y^i}{y^j}(y+\epsilon w)
\]
for some $\epsilon$ with $0<\epsilon<1$. The final term is
non-negative, so
\[
F(y+w)\geq F(y)+w^i\fpd{F}{y^i}(y),
\]
and similarly
\[
F(y-w)\geq F(y)-w^i\fpd{F}{y^i}(y),
\]
where we must now assume that $0$ does not lie in the segment
$[y-w,y+w]$. Summing,
\[
F(y+w)+F(y-w)\geq 2F(y).
\]
On setting $y_1=\onehalf(y+w)$ and $y_2=\onehalf(y-w)$ we obtain
the triangle inequality, for all pairs of
points $y_1,y_2$ such that the segment $[y_1,y_2]$ does not
contain the origin. Now by continuity $F(0)=0$, so the triangle
inequality holds self-evidently if either argument is $0$.
Suppose that $0$ is an interior point of the segment
$[y_1,y_2]$. Consider first the case $y_1+y_2=0$. For any $y\neq
0$ and any positive $s$, $F(sy+y_1)+F(sy+y_2)\geq 2sF(y)$, and
on letting $s$ tend to 0 we obtain $F(y_1)+F(y_2)\geq0$. Now
suppose that $y_2=-ky_1$ with $k>0$, and furthermore suppose
that $k<1$ (otherwise reverse the roles of $y_1$ and $y_2$).
Then $y_1+y_2=(1-k)y_1$, and
\[
F(y_1)+F(y_2)=F(y_1)+kF(-y_1)\geq (1-k)F(y_1)=F(y_1+y_2).
\]
Thus the triangle inequality holds for all $y_1,y_2$.

For the fundamental inequality replace $y+w$ by $z$ to obtain
\[
F(z)\geq F(y)+(z^i-y^i)\fpd{F}{y^i}(y)=z^i\fpd{F}{y^i}(y).
\]
We have to exclude $y=0$ since $F$ won't be differentiable
there. We also have to assume that the segment $[y,z]$ does not
include the origin, that is, that $z$ is not a negative multiple
of $y$. But if $z=-ky$ with $k>0$ then
\[
F(z)=kF(-y)\geq -kF(y)=z^i\fpd{F}{y^i}(y).
\qedhere
\]
\end{proof}

Note that  by positive homogeneity $F$ can only
be constant if it is identically zero.

\begin{cor}
If $F$ is not constant then it must take positive values.
\end{cor}

\begin{proof}
If $F$ is not constant there must be a point $y$ at which
the gradient of $F$ is non-zero. Then for all $z$ in an open
half-space, $z^i\partial F/\partial y^i(y)> 0$; and $F(z)>0$ for
all such $z$.
\end{proof}

We assume for the rest of this section that the Hessian of $F$
is positive quasi-definite.

\begin{lem}\label{constant}
It cannot be the case that for fixed non-zero $z$,
\[
z^i\fpd{F}{y^i}(y)
\]
is constant (i.e.\ takes the same value for all $y\neq 0$).
\end{lem}

\begin{proof}
If it were constant it would  follow that
\[
z^j\spd{F}{y^i}{y^j}(y)=0
\]
for all $y$, which cannot be the case when the Hessian is
positive quasi-definite if $z\neq 0$.
\end{proof}

\begin{lem}\label{strict}
Except when $z$ is a positive scalar multiple of $y$, the fundamental
inequality is strict:
\[
 F(z)>z^i\fpd{F}{y^i}(y).
\]
\end{lem}

\begin{proof}
Evidently equality holds when $z$ is a positive scalar multiple
of $y$. Moreover, it follows from the second mean value theorem
that the inequality is strict otherwise, except possibly when
$z$ is a negative scalar multiple of $y$. From the proof of
Proposition \ref{tri} we know that if $z=-ky$ with $k>0$ then
\[
F(z)=kF(-y)\geq -kF(y)=z^i\fpd{F}{y^i}(y).
\]
So the question is whether we can have
$F(-y)=-F(y)$ with $y\neq 0$. Now for every $w\neq 0$,
\[
F(-y)\geq -y^i\fpd{F}{y^i}(w),\quad
-F(y)\leq -y^i\fpd{F}{y^i}(w);
\]
so if $F(-y)=-F(y)=-c$ say then
\[
y^i\fpd{F}{y^i}(w)=c\mbox{ for all $w\neq 0$.}
\]
But by Lemma \ref{constant} this cannot hold. Thus $F(-y)>-F(y)$
and the inequality is strict in this case too.
\end{proof}

\begin{prop}\label{abshom}
If $F$ is absolutely homogeneous and its Hessian is positive
quasi-definite then it is everywhere positive on $\R^n-\{0\}$.
\end{prop}

\begin{proof}
It follows from the previous proof that if $F(y)=0$ (where
$y\neq 0$) then $F(-y)>0$. But if $F$ is absolutely homogeneous
and $F(y)=0$ then it would have to be the case that $F(-y)=0$
also. So $F$ can never vanish. We know that $F$ must take
positive values. So it must be positive everywhere on
$\R^n-\{0\}$.
\end{proof}
\noindent This result has appeared previously in \cite{hil4}.

\begin{prop}
If $F$ is positively homogeneous and its Hessian is positive
quasi-definite then there is a linear function $y\mapsto
\alpha_iy^i$ such that if $\tilde{F}(y)=F(y)+\alpha_iy^i$ then
$\tilde{F}$ is everywhere positive on $\R^n-\{0\}$.
\end{prop}

\begin{proof}
Fix $z\in\R^n-\{0\}$ and define $\bar{F}$ by
\[
\bar{F}(y)=F(y)-y^i\fpd{F}{y^i}(z).
\]
Then by Lemma \ref{strict}, $\bar{F}>0$ on $\R^n-\{0\}$, except
along the ray through $z$; $\bar{F}(z)=0$. Now consider the
restriction of $\bar{F}$ to the unit coordinate Euclidean sphere
$S$. Denote by $H_z$ the closed hemisphere opposite the point
where the ray through $z$ intersects $S$, that is, the southern
hemisphere where $\hat{z}=z/|z|$ is the north pole. Then
$\bar{F}$ is positive on $H_z$, which is compact, so there is
some $k>0$ such that $\bar{F}(y)\geq k$ for $y\in H_z$. Now set
$\varphi(y)=\onehalf k\hat{z}_iy^i$ where
$\hat{z}_i=\delta_{ij}\hat{z}^j$. Note that $-\onehalf
k\leq\varphi(y)\leq\onehalf k$ for $y\in S$:\ in fact $\varphi$
takes its minimum value $-\onehalf k$ at the south pole, is zero
on the equator, and is positive on the open northern hemisphere
$S-H_z$. Set $\tilde{F}=\bar{F}+\varphi$. Then $\tilde{F}$ is
positive on $S-H_z$, since $\bar{F}$ is non-negative and
$\varphi$ is positive; and also on $H_z$, since $\bar{F}\geq k$
and $\varphi\geq -\onehalf k$. Thus $\tilde{F}$ is positive
everywhere on $S$, and so everywhere on $\R^n-\{0\}$; and
$\tilde{F}(y)=F(y)+\alpha_iy^i$ where
\[
\alpha_i=\onehalf k\hat{z}_i-\fpd{F}{y^i}(z)
\]
for the chosen non-zero $z$.
\end{proof}

\begin{thm}\label{totderiv}
Let $F$ be a positively homogeneous function on $\TMO$ such that
for every $x\in M$ the Hessian of $F$ on $\TxMO$ is positive
quasi-definite. Then for any $x_0\in M$ there is a neighbourhood
$U$ of $x_0$ in $M$ and a function $\tilde{F}$ defined on
$T^\circ U$ such that $\tilde{F}$ is a Finsler function which
differs from $F$ by a total derivative. If $F$ is absolutely,
rather than just positively, homogeneous then it itself is a
Finsler function.
\end{thm}

\begin{proof}
Take any point $x_0\in M$ and choose coordinates at $x_0$. By
the previous proposition there is a (constant) covector
$\alpha$ such that if $\tilde{F}(x,y)=F(x,y)+\alpha_iy^i$ then
$\tilde{F}(x_0,y)>0$ for $y\neq 0$. It follows that there is an
open neighbourhood $U$ of $x_0$ within the coordinate patch such
that $\tilde{F}(x,y)>0$ for all $x\in U$ and $y\in\TxMO$. (To
see this, suppose the contrary. Consider $\tilde{F}$ as function
on $P\times S$ where $P$ is the coordinate patch and $S$, as
before, is the unit Euclidean sphere. Then there is a sequence
of points $\{(x_r,y_r)\}$ in $P\times S$ with $x_r\to x_0$ such
that $\tilde{F}(x_r,y_r)\leq 0$ for all $r$. But $\{y_r\}$,
being a sequence in a compact set, has a convergent subsequence,
which converges to $y_0$ say. Then $\tilde{F}(x_0,y_0)\leq 0$,
which is a contradiction.) So $\tilde{F}$ is a Finsler function
on $T^\circ U$, and it differs from $F$ by $\alpha_iy^i$ which
is the total derivative of the local function $\alpha_ix^i$.

If $F$ is absolutely homogeneous then it is everywhere positive
by Proposition \ref{abshom}, and therefore a Finsler function.
\end{proof}

\section{The Euler-Lagrange equations and the Rapcs\'{a}k
conditions}

For a function $F$ on $\TMO$ which is positively homogeneous of
degree 1, the Euler-Lagrange equations
\[
\Gamma\left(\fpd{F}{y^i}\right)-\fpd{F}{x^i}=0
\]
have the property that if there is a solution $\Gamma$ which is a
spray, then every spray projectively equivalent to $\Gamma$ is also a
solution.

At a basic level therefore solving the projective metrizability
problem for a given projective class of sprays involves finding
a positively-homogeneous function $F$ such that the
Euler-Lagrange equations are satisfied for any spray $\Gamma$ in
the class. These equations, considered as conditions on $F$,
constitute (in the terminology of \cite{Szi,SzVa}) the first set
of {\em Rapcs\'{a}k conditions}. There is a second equivalent
set which is given in the next proposition. (As a matter of fact
there are many versions of the Rapcs\'{a}k conditions to be
found in the literature:\ \cite{Bacso} contains seven. The two
we give are the ones most useful for our purposes.)

\begin{prop}
The positively-homogeneous function $F$ satisfies the first
Rapcs\'{a}k conditions if and only if
\[
\spd{F}{x^j}{y^i}-\Gamma^k_j\spd{F}{y^i}{y^k}=
\spd{F}{x^i}{y^j}-\Gamma^k_i\spd{F}{y^j}{y^k}.
\]
\end{prop}

\begin{proof}
The first Rapcs\'{a}k conditions are
\[
y^k\spd{F}{x^k}{y^i}-2\Gamma^k\spd{F}{y^i}{y^k}-\fpd{F}{x^i}=0.
\]
On differentiating with respect to $y^j$ one obtains
\[
\Gamma\left(\spd{F}{y^i}{y^j}\right)+\spd{F}{x^j}{y^i}
-2\Gamma^k_j\fpd{F}{y^i}{y^k}-\spd{F}{x^i}{y^j}=0.
\]
The part of this which is skew in $i$ and $j$ gives the required
conditions. Conversely, if the given conditions hold and $F$ is
homogeneous then contraction with $y^j$ produces the first
Rapcs\'{a}k conditions.
\end{proof}
\noindent It is of interest, though not of any significance here, that the
second Rapcs\'{a}k conditions follow from the first regardless
of whether or not $F$ is homogeneous; homogeneity is required
for the converse step only.

The second Rapcs\'{a}k conditions may be written
\[
H_j(\theta_i)=H_i(\theta_j)\quad\mbox{where }\theta_i=\fpd{F}{y^i}.
\]
The symmetric part of the equation whose skew part gives the
second Rapcs\'{a}k conditions may be written
\[
(\nabla h)_{ij}=0\quad \mbox{where }h_{ij}=\spd{F}{y^i}{y^j},
\]
where $\nabla$ is the dynamical covariant derivative operator
associated with $\Gamma$. These conditions have important roles
to play in the next section.

\section{The Helmholtz-like conditions for a multiplier}\label{Hholtz}

This section is concerned with conditions on a multiplier
necessary and sufficient for it to be the Hessian of a Finsler
function for a given projective class of sprays.

Initially we work in a coordinate neighbourhood $U$ of $M$.
Without essential loss of generality we may and shall assume
that the Poincar\'{e} Lemma holds on $U$:\ in fact we shall
assume that $U$ is contractible. We further assume, for reasons
that will shortly become clear, that $\dim M\geq 3$; and indeed
we make this assumption for most of the rest of this paper,
until the discussion of the two-dimensional case in Section 7 in
fact.

We denote a putative Hessian of a Finsler function by $h$, and
its components by $h_{ij}$, to distinguish it from the putative
Hessian of an energy function, which is written $g$ or $g_{ij}$
in the conventional way. We shall say that $h$ is {\em quasi-regular\/}
if $h_{ij}v^j=0$ if and only if $v^i=ky^i$ for some constant $k$.
We shall call a positively-homogeneous function $F$ whose Hessian is
quasi-regular a {\em pseudo-Finsler function}.

\begin{lem}\label{first}
If $h_{ij}$ is defined and smooth on $T^\circ U$, is symmetric
and satisfies $h_{ij}y^j=0$ and
\[
\fpd{h_{ij}}{y^k}=\fpd{h_{ik}}{y^j}
\]
then
\begin{enumerate}
\item $h_{ij}$ is positively homogeneous of degree $-1$;
\item there is a smooth function $F$ on $T^\circ U$, positively
homogeneous of degree 1, such that
\[
h_{ij}=\spd{F}{y^i}{y^j}.
\]
\end{enumerate}
\end{lem}

\begin{proof}
1. We have
\[
0=\fpd{}{y^j}(h_{ik}y^k)=\fpd{h_{ik}}{y^j}y^k+h_{ij}=
\fpd{h_{ij}}{y^k}y^k+h_{ij}
\]
as required.\\[5pt]
2. For each $i$ we may consider $h_{ij}dy^j$ as a 1-form on
$\R^n-\{0\}$, $n=\dim M$, with coordinates $y^i$, where we
regard the $x^i$ as parameters. The condition
\[
\fpd{h_{ij}}{y^k}=\fpd{h_{ik}}{y^j}
\]
says that this 1-form is closed. For $\dim M\geq 3$ it follows
that $h_{ij}dy^j$ is exact, that is, there are functions
$\bar{F}_i$ depending smoothly on $y^i$ for $y\neq 0$, and
smoothly also on the parameters $x^i$, such that
\[
h_{ij}=\fpd{\bar{F}_i}{y^j}.
\]
But by the symmetry of $h_{ij}$, $\bar{F}_i dy^i$ is also
closed, and hence exact. By a similar argument there is a smooth
function $\bar{F}$ on $T^\circ U$ such that
\[
h_{ij}=\spd{\bar{F}}{y^i}{y^j}.
\]
Then
\[
\fpd{}{y^i}\left(y^j\fpd{\bar{F}}{y^j}-\bar{F}\right)=h_{ij}y^j=0,
\]
so that $y^i\partial \bar{F}/\partial y^i-\bar{F}$ is a function of
$x^i$ alone, say $f$. Then $F=\bar{F}+f$ has the required properties;
and so indeed  has any function differing from $F$ by a function
linear in the $y^i$.
\end{proof}
\noindent The first result also follows from the second, of course, but the
direct proof is available if one does not want to invoke the
existence of $F$.

Note that we require that $\dim M=n\geq 3$ to be able to
conclude that a closed 1-form on $\R^n-\{0\}$ is exact, and
hence that $F(x,y)$ is defined for all $y\neq 0$. It is for this
reason that we make this assumption for most of the rest of the
paper.

\begin{lem}
Let $\Gamma$, $\tilde{\Gamma}$ be projectively-equivalent sprays,
$\nabla$, $\tilde{\nabla}$ be the corresponding dynamical covariant
derivative operators. If $h_{ij}$ is symmetric and homogeneous of
degree $-1$ then $(\tilde{\nabla}h)_{ij}=(\nabla h)_{ij}$.
\end{lem}

\begin{proof}
With $\tilde{\Gamma}=\Gamma-2P\Delta$ and $P_i=V_i(P)$ we have
\begin{align*}
(\tilde{\nabla}h)_{ij}
&=\tilde{\Gamma}(h_{ij})-\tilde{\Gamma}_i^kh_{kj}-\tilde{\Gamma}_j^kh_{ik}\\
&=\Gamma(h_{ij})-2P\Delta(h_{ij})\\
&\quad\mbox{}-(\Gamma_i^k+P\delta_i^k+P_iy^k)h_{kj}
-(\Gamma_j^k+P\delta_j^k+P_jy^k)h_{ik}\\
&=(\nabla h)_{ij}+2Ph_{ij}-Ph_{ij}-Ph_{ij}\\
&=(\nabla h)_{ij}.\qedhere
\end{align*}
\end{proof}
\noindent Thus $\nabla h$ is projectively invariant for a given projective
class.

The following lemma involves the curvature tensors $R_j^i$ and
$R^i_{jk}$ of  a spray $\Gamma$.

\begin{lem}\label{hR}
If $h_{ij}$ satisfies the conditions of Lemma~\ref{first} then the
following conditions are equivalent to each other:
\begin{align*}
&h_{ik}R^k_j=h_{jk}R^k_i\\
&h_{il}R^l_{jk}+h_{jl}R^l_{ki}+h_{kl}R^l_{ij}=0.
\end{align*}
\end{lem}

\begin{proof}
To obtain the second expression, differentiate
$h_{il}R^l_j-h_{jl}R^l_i$ with respect to $y^k$ to obtain
\[
\fpd{h_{il}}{y^k}R^l_j+h_{il}\fpd{R^l_j}{y^k}
-\fpd{h_{jl}}{y^k}R^l_i+h_{jl}\fpd{R^l_i}{y^k}=0,
\]
and add the two similar expressions obtained by cyclically permuting
$i$, $j$ and $k$.  Terms involving derivatives of the $h_{ij}$ cancel
in pairs due to the symmetry condition they satisfy, while pairs of
terms involving derivatives of the $R^l_i$ give the $R^l_{ij}$.

To obtain the first from the second, simply contract with $y^k$.
\end{proof}
\noindent This lemma is also to be found in \cite{Bacso}.

The Weyl curvature $W^i_j$ of a spray $\Gamma$
is defined by
\[
W^i_j=R^i_j-R\delta^i_j-\rho_jy^i\quad\mbox{where }
R=\frac{1}{n-1}R^k_k \mbox{ and }
\rho_j=\frac{1}{n+1}\left(\fpd{R^k_j}{y^k}-\fpd{R}{y^j}\right).
\]
The Weyl curvature is projectively invariant. See \cite{Shen}
for details.

\begin{lem}\label{Weyl}
If $h_{ij}$ satisfies the conditions of Lemma~\ref{first} then the
condition $h_{ik}R^k_j=h_{jk}R^k_i$ is equivalent to
$h_{ik}W^k_j=h_{jk}W^k_i$, and in particular is projectively invariant.
\end{lem}

\begin{proof}
We have
\begin{align*}
h_{ik}W^k_j-h_{jk}W^k_i&=
h_{ik}(R_j^k-R\delta_j^k-\rho_jy^k)-h_{jk}(R_i^k-R\delta_i^k-\rho_iy^k)\\
&=h_{ik}R_j^k-Rh_{ij}-h_{jk}R_i^k+Rh_{ij}\\
&=h_{ik}R^k_j-h_{jk}R^k_i
\end{align*}
as claimed.
\end{proof}

\begin{lem}\label{rap2}
Let $\bar{F}$ be a smooth function on
$T^\circ U$ which is positively homogeneous and satisfies
\[
(\nabla h)_{ij}=0\quad \mbox{ and }\quad
h_{ik}W^k_j=h_{jk}W^k_i
\]
where
\[
h_{ij}=\spd{\bar{F}}{y^i}{y^j}.
\]
Then there is a smooth positively-homogeneous function $F$ on
$T^\circ U$, with the same Hessian as $\bar{F}$, which satisfies
the second Rapcs\'{a}k conditions.
\end{lem}

\begin{proof}
We show first that $H_i(\bar{\theta}_j)-H_j(\bar{\theta}_i)$ is
independent of the $y^k$, where
$\bar{\theta}_i=\partial\bar{F}/\partial y^i$. Now
\[
\vf{y^k}\big(H_i(\bar{\theta}_j)\big)=
H_i(h_{jk})-\conn lik h_{jl},\quad
\conn kij=\fpd{\Gamma^k_i}{y^j}=\conn kji.
\]
It is a simple and well-known consequence of the first
assumption, together with the evident fact that
\[
\fpd{h_{ij}}{y^k}=\fpd{h_{ik}}{y^j},
\]
that
\[
H_i(h_{jk})-\conn likh_{jl}=
H_j(h_{ik})-\conn ljkh_{il},
\]
whence $H_i(\bar{\theta}_j)-H_j(\bar{\theta}_i)$ is independent of the
$y^k$. Thus
\[
\big(H_i(\bar{\theta}_j)-H_j(\bar{\theta}_i)\big)dx^i\wedge dx^j
\]
is a basic 2-form, say $\chi$. We show next that $\chi$ is closed. In
computing $d\chi$ we may replace partial derivatives with respect to
$x^k$ with $H_k$. We have
\[
\oplus H_k\big(H_i(\bar{\theta}_j)-H_j(\bar{\theta}_i)\big)
=\oplus[H_j,H_k](\bar{\theta}_i)=-\oplus R^l_{jk}h_{il}
\]
where $\oplus$ indicates the cyclic sum over $i$, $j$ and $k$.
But this vanishes if $h_{ik}W^k_j=h_{jk}W^k_i$ by Lemmas
\ref{hR} and \ref{Weyl}. So $\chi$ is closed, and hence exact.
Choose $\psi=\psi_idx^i$ such that $\chi=d\psi$, and set
$F=\bar{F}-\psi_iy^i$, so that $\theta_i=\bar{\theta}_i-\psi_i$.
Then
\[
\big(H_i(\theta_j)-H_j(\theta_i)\big)dx^i\wedge dx^j=\chi-d\psi=0,
\]
and $F$ satisfies the second Rapcs\'{a}k conditions.
\end{proof}

After all these preliminaries we can now give the main result of
this section.

\begin{thm}\label{mult}
Let $U\subset M$ be any contractible open subset of a coordinate
patch. Given a projective class of sprays, the following
conditions are necessary and sufficient for the existence of a
positively-homogeneous function $F$ on $T^\circ U$, such that
every spray in the class satisfies the Euler-Lagrange equations
for $F$:\ there is a tensor $h$ whose components satisfy
\begin{align*}
h_{ji}&=h_{ij}\\
h_{ij}y^j&=0\\
\fpd{h_{ij}}{y^k}&=\fpd{h_{ik}}{y^j}\\
(\nabla h)_{ij}&=0\\
h_{ik}W^k_j&=h_{jk}W^k_i,
\end{align*}
where $\nabla$ is the dynamical covariant derivative operator of any
spray in the class.
\end{thm}

\begin{proof}
If such a function $F$ exists then its Hessian satisfies these
conditions. Conversely, from the first three conditions and
Lemma \ref{first} there is a smooth positively-homogeneous
function $\bar{F}$ on $T^\circ U$ whose Hessian is $h$. Take any
$\Gamma$ in the projective class. Then by Lemma \ref{rap2} and
the remaining conditions there is a smooth positively
homogeneous function $F$ on $T^\circ U$, whose Hessian is also
$h$, which satisfies the second Rapcs\'{a}k conditions. By the
equivalence of the first and second Rapcs\'{a}k conditions
$\Gamma$, and hence every spray in the class, satisfies the
Euler-Lagrange equations for $F$.
\end{proof}

\begin{cor}
Let $U\subset M$ be any contractible open subset of a coordinate patch.
If there is a tensor $h$ which satisfies the conditions of
Theorem~\ref{mult} and $F$ is a corresponding
positively-homogeneous function on $T^\circ U$ such that every spray in
the class satisfies the Euler-Lagrange equations for $F$ then
another function $\tilde{F}$ has the same properties if and
only if differs from $F$ by the total derivative of a function on $U$.
\end{cor}

\begin{proof}
The function $\tilde{F}$ has the same Hessian as $F$ and is
positively homogeneous, so $\tilde{F}=F+\alpha_i y^i$ for some
functions $\alpha_i$ on $U$. If $\Gamma$ satisfies the
Euler-Lagrange equations for both $F$ and $\tilde{F}$ we must
have
\[
y^j\fpd{\alpha_i}{x^j}-y^j\fpd{\alpha_j}{x^i}=0,
\]
or equivalently $\alpha_i dx^i$ is closed and therefore exact; and if
$\alpha_idx^i=d\phi$ then $\tilde{F}-F$ is the total derivative of $\phi$.
\end{proof}

\section{Some global results}\label{glob}

We now consider the problem of extending these results from
coordinate neighbourhoods in $M$ to the whole of $M$. We shall
work with an open covering
$\mathfrak{U}=\{U_\lambda:\lambda\in\Lambda\}$ of $M$ by
coordinate patches. We shall assume that $\mathfrak{U}$ has the
property that every $U_\lambda$, and every non-empty
intersection of finitely many of the $U_\lambda$, is
contractible. A covering with this property is known as a good
covering; it can be shown (see Appendix B) that every manifold
over which is defined a spray admits good open coverings by
coordinate patches.

Let us assume that there is a tensor field $h$  which
satisfies the conditions of Theorem~\ref{mult}, and in addition
is everywhere quasi-regular. Then for each $U_\lambda$ there is a
pseudo-Finsler function $F_\lambda$ defined on $T^\circ
U_\lambda$; and on $U_\lambda \cap U_\mu$, if it is non-empty,
there is a function $\phi_{\lambda\mu}$, determined up to the
addition of a constant, such that
\[
F_\lambda-F_\mu= y^i\fpd{\phi_{\lambda\mu}}{x^i}.
\]
On $U_\lambda \cap U_\mu \cap U_\nu$, if it is non-empty,
$\phi_{\mu\nu}-\phi_{\lambda\nu}+\phi_{\lambda\mu}$ is a
constant, say $k_{\lambda\mu\nu}$ (add the expressions for
$F_\lambda-F_\mu$ etc., and note that $U_\lambda \cap U_\mu \cap
U_\nu$, being contractible, is connected). Furthermore, for any
four members $U_\kappa,U_\lambda,U_\mu,U_\nu$ of $\mathfrak{U}$
whose intersections in threes are non-empty
\[
 k_{\lambda\mu\nu}-k_{\kappa\mu\nu}+k_{\kappa\lambda\nu}
-k_{\kappa\lambda\mu}=0.
\]
That is to say, $k$ satisfies a cocycle condition. If $k$ is in
fact a coboundary we can modify each $\phi_{\lambda\mu}$ by the
addition of a constant, so that (after modification)
$\phi_{\mu\nu}-\phi_{\lambda\nu}+\phi_{\lambda\mu}=0$.

We now show that if, for an open covering $\mathfrak{U}$ of
$M$ by coordinate patches we can find functions
$\phi_{\lambda\nu}$ which satisfy this cocycle condition
then there is a globally-defined pseudo-Finsler function $F$ on
$\TMO$ for the given projective class of sprays.

Since by assumption $M$ is paracompact the covering
$\mathfrak{U}$ of $M$ admits a locally finite refinement
$\mathfrak{V}=\{V_\alpha:\alpha\in A\}$. (The covering
$\mathfrak{V}$ may not be good, but this will not matter for
what follows.) There is a partition of unity subordinate to
$\mathfrak{V}$, that is, for every $\alpha$ there is a smooth function
$f_\alpha$ such that $\supp(f_\alpha)\subset V_\alpha$, $0\leq
f_\alpha\leq 1$ and $\sum_\alpha f_\alpha=1$. By local
finiteness, for every $x\in M$, $\sum_\alpha f_\alpha(x)$ is a
finite sum. Indeed, since $\supp(f_\alpha)\subset V_\alpha$
\[
\sum_\alpha f_\alpha(x)
=\sum_{\alpha:x\in V_\alpha}f_\alpha(x).
\]

Since $\mathfrak{V}$ is a refinement of $\mathfrak{U}$ we can
define a map $\sigma: A\to\Lambda$ such that $V_\alpha\subseteq
U_{\sigma(\alpha)}$. For every $\alpha,\beta\in A$ such that
$V_\alpha\cap V_\beta$ is non-empty we define a real-valued
function $\phi_{\alpha\beta}$ by
\[
\phi_{\alpha\beta}=
\phi_{\sigma(\alpha)\sigma(\beta)}|_{V_\alpha\cap V_\beta}.
\]
Then the $\phi_{\alpha\beta}$ satisfy the cocycle condition, namely
\[
\phi_{\beta\gamma}-\phi_{\alpha\gamma}+\phi_{\alpha\beta}=0,
\]
since the $\phi_{\lambda\mu}$ do. Moreover, if we set
$F_\alpha=F_{\sigma(\alpha)}|_{\tau^{-1}V_\alpha}$ then
\[
F_\alpha-F_\beta=y^i\fpd{\phi_{\alpha\beta}}{x^i}.
\]

\begin{prop}
Under these assumptions there is a globally-defined function $F$
on $\TMO$ such that $F_\alpha$ differs from
$F|_{\tau^{-1}V_\alpha}$ by a total derivative.
\end{prop}

\begin{proof}
Let $V_\alpha$ be any member of the covering $\mathfrak{V}$. For
$x\in V_\alpha$ set
\[
\psi_\alpha(x)=
\sum_{\gamma:x\in V_\gamma}f_\gamma(x)\phi_{\alpha\gamma}(x).
\]
Consider $\psi_\alpha(x)-\psi_\beta(x)$ where
$x\in V_\alpha\cap V_\beta$. We have
\begin{align*}
\psi_\alpha(x)-\psi_\beta(x)&=
\sum_{\gamma:x\in V_\gamma}f_\gamma(x)(\phi_{\alpha\gamma}(x)
-\phi_{\beta\gamma}(x))\\
&=\sum_{\gamma:x\in V_\gamma}f_\gamma(x)\phi_{\alpha\beta}(x)\\
&=\left(\sum_{\gamma:x\in V_\gamma}f_\gamma(x)\right)
\phi_{\alpha\beta}(x)\\
&=\phi_{\alpha\beta}(x).
\end{align*}
Then for any two members of the covering $\mathfrak{V}$ with
non-empty intersection
\[
F_\alpha-y^i\fpd{\psi_\alpha}{x^i}=F_\beta-y^i\fpd{\psi_\beta}{x^i}
\]
on $V_\alpha\cap V_\beta$; and we can define $F$ consistently by
\[
F(x,y)=F_\alpha(x,y)-y^i\fpd{\psi_\alpha}{x^i}(x)
\]
where $V_\alpha$ is any member of the covering containing $x$.
\end{proof}
\noindent The first part of this proof is a simple special case of the
argument used to show that \v{C}ech cohomology is a sheaf
cohomology theory in \cite{Warner}.

So for $\dim M \geq 3$ the multiplier argument gives a global
pseudo-Finsler function provided the cocycle $k$ is a
coboundary. Now $k$ is in fact an element of the \v{C}ech
cochain complex for the covering $\mathfrak{U}$ with values in
the constant sheaf $M\times\R$. This cocycle will certainly be a
coboundary if the corresponding cohomology group on $M$, namely
$\check{H}^2(\mathfrak{U},\mathcal{R})$, is zero (we use the
notation of \cite{Warner}). It follows from the fact that
$\mathfrak{U}$ is a good open covering that this cohomology
group is isomorphic to the de Rham cohomology group $H^2(M)$. We
have the following result.

\begin{thm}
If $F$ is a (global) Finsler function on $\TMO$ then its Hessian
$h$ satisfies the conditions of Theorem~\ref{mult} for the
sprays of its geodesic class, and is in addition positive
quasi-definite. Conversely, suppose given a
projective class of sprays on $\TMO$. If there is a
tensor field $h$ which everywhere satisfies the
conditions of Theorem~\ref{mult} and is in addition positive
quasi-definite, and if $H^2(M)=0$, then the projective class is
the geodesic class of a global pseudo-Finsler function, and of a
local Finsler function over a neighbourhood of any point of $M$.
\end{thm}

We give examples of sprays which admit a global
pseudo-Finsler function but only local Finsler functions in
Section 8, so the result above would appear to be the best possible.

The situation in the case of a reversible spray, or without loss
of generality in the light of Proposition \ref{rev} a strictly
reversible spray, is much more clear cut. We observe first that
if a spray $\Gamma$ satisfies the Euler-Lagrange equations of a
positively-homogeneous function $F$, and we set
$\bar{F}=\rho^*F$ (where $\rho$ is the reflection map), then
$\bar{\Gamma}$, the reverse of $\Gamma$, satisfies the
Euler-Lagrange equations of the (positively-homogeneous)
function $\bar{F}$. So if $\Gamma$ is strictly reversible, it
satisfies the Euler-Lagrange equations for both $F$ and
$\bar{F}$, and therefore for their sum. But if $F$ and $\bar{F}$
are positively homogeneous, $F+\bar{F}$ is absolutely
homogeneous: for if $k<0$
\begin{align*}
F(ky)+\bar{F}(ky)&=F(-|k|y)+\bar{F}(-|k|y)\\
&=|k|(F(-y)+\bar{F}(-y))=|k|(\bar{F}(y)+F(y)).
\end{align*}
So for a strictly reversible spray, on $U_\lambda$ there is a
pseudo-Finsler function $F_\lambda$ which is absolutely
homogeneous. Now two functions (such as $F_\lambda$, $F_\mu$)
which are known to differ by a total derivative and are both
absolutely homogeneous must be equal. In the light of these
remarks and the final assertion of Theorem \ref{totderiv} we
have the following result.

\begin{thm}\label{abs}
The projective class of a reversible spray on $\TMO$
is the geodesic class of a globally-defined absolutely-homogeneous
Finsler function if and only if there is a tensor
$h$ which satisfies the conditions of Theorem~\ref{mult} and is
in addition positive quasi-definite.
\end{thm}

\section{The two-dimensional case}

It is well known that in two dimensions every spray is locally
projectively metrizable (see for example \cite{APB,BM}). `Locally'
here, however, refers to the fibre as well as to the base:\ the
claim is only that a y-local Finsler function exists. We may
seek to establish the existence of a Finsler, or at least
pseudo-Finsler, function as follows.

\begin{lem}
Let $h$ be a symmetric smooth  tensor field on $\R^2-\{0\}$ which is positively homogeneous of degree $-1$ and satisfies $h_{ij}y^j=0$. Then
\[
\fpd{h_{ij}}{y^k}=\fpd{h_{ik}}{y^j}.
\]
\end{lem}

\begin{proof}
We have
\begin{align*}
h_{11}y^1+h_{12}y^2&=0\\
h_{21}y^1+h_{22}y^2&=0.
\end{align*}
Thus by symmetry $h$ has just one independent component, which
we may take to be $h_{12}$; this must be non-zero except where
$y^1=0$ or $y^2=0$ (separately) for $h$ to be non-trivial.
Moreover $h_{12}$ must satisfy the homogeneity condition
\[
y^1\fpd{h_{12}}{y^1}+y^2\fpd{h_{12}}{y^2}+h_{12}=0.
\]
On differentiating the equation $h_{11}y^1+h_{12}y^2=0$ with respect to
$y^2$ we obtain
\[
y^1\fpd{h_{11}}{y^2}+y^2\fpd{h_{12}}{y^2}+h_{12}=0,
\]
whence by homogeneity of $h_{12}$
\[
y^1\left(\fpd{h_{11}}{y^2}-\fpd{h_{12}}{y^1}\right)=0.
\]
Thus
\[
\fpd{h_{11}}{y^2}=\fpd{h_{12}}{y^1}
\]
except possible where $y^1=0$. But by continuity this must hold
where $y^1=0$ (but $y^2\neq 0$) also. This is the first of the
two non-trivial cases of the equations it is required to prove;
the other is proved similarly.
\end{proof}

\begin{lem}
Let $U$ be a coordinate neighbourhood in a two-dimensional
manifold on which is defined a spray $\Gamma$. Suppose that a
symmetric tensor $h_{ij}$ satisfies $(\nabla h)_{ij}=0$. If
$\eta_i=h_{ij}y^j$ then $\nabla \eta=0$, and if
$\lambda_{ij}=\Delta(h_{ij})+h_{ij}$ then $\nabla\lambda=0$.
\end{lem}

\begin{proof}
The first follows immediately from the fact that $\nabla y=0$
(that is to say, the dynamical covariant derivative of the total
derivative $y^i\partial/\partial x^i$ vanishes). For the second
we use the facts that $[\Delta,\Gamma]=\Gamma$ and that
$\Delta(\Gamma^i_j)=\Gamma^i_j$ to show that
$[\Delta,\nabla]=\nabla$.
\end{proof}

The Weyl tensor $W^i_j$ vanishes when $\dim M=2$. In order for
the Helmholtz conditions as given in Theorem~\ref{mult} to be
satisfied it is therefore enough that $\nabla h=0$ everywhere,
and that $\eta=0$ and $\lambda=0$ on some cross-section of the
flow of $\Gamma$, which we must assume to contain, for every $y$
in it, the ray through $y$, for the latter condition to make
sense. Indeed, we can always find a multiplier $h$ locally in
$\TMO$ by specifying its values on such a cross-section of the
flow of $\Gamma$, arbitrarily subject to the conditions that
$\eta=0$ and $\lambda=0$, and requiring that $\nabla h=0$:\ this
is a first-order differential equation along each integral curve
of $\Gamma$, and so the value of $h$ along the curve is
determined by the specified initial value. Moreover, $\eta$ and
$\lambda$ satisfy similar first-order differential equations
along the integral curves of $\Gamma$; and now we can use the
fact that the equations are linear to show that since $\eta=0$
initially, $\eta=0$ everywhere, and likewise for $\lambda$.
Since $\Gamma$ never vanishes on $\TMO$ we can always find a
local cross-section to its flow. Indeed, we can take a local
two-dimensional submanifold of $\TMO$ which is transverse to
$\Gamma$, and extend it to a local cross-section by including,
for every $y$ in it, the ray through $y$. This, together with
the evident y-local version of Theorem~\ref{mult}, establishes
the result stated above.

\begin{thm}
When $\dim M=2$, every spray is locally projectively metrizable.
\end{thm}

It is moreover easy to specify the freedom in the choice of multiplier.

\begin{prop}\label{2d}
Suppose that in two dimensions $h$ satisfies the Helmholtz
conditions. Let $f$ be a function on $\TMO$ such that $Z(f)=0$
for every vector field $Z\in\langle\Delta,\Gamma\rangle$ (so
that $f$ is homogeneous of degree 0, and constant along the
integral curves of any spray in the projective class). Then
$\hat{h}=fh$ also satisfies the Helmholtz conditions.
Conversely, any pair of tensors which satisfy the Helmholtz
conditions are so related.
\end{prop}

\begin{proof}
It is clear that $\hat{h}$ satisfies the algebraic conditions
and is homogeneous of degree $-1$. We have
\[
(\nabla\hat{h})_{ij}=f(\nabla h)_{ij}+\Gamma(f)h_{ij}=0.
\]
Conversely, if $h_{ij}y^j=0=\hat{h}_{ij}y^j$ then $\hat{h}$ must be a
scalar multiple of $h$, say $\hat{h}=fh$. Then by
homogeneity we must have $\Delta(f)=0$, and if
$\nabla\hat{h}=0=\nabla h$ then $\Gamma(f)=0$ also.
\end{proof}

The results above are y-local. But even if there is a y-global
multiplier this by itself is not enough to guarantee the
existence of a y-global positively-homogeneous function $F$ on
$\TMO$ of which it is the Hessian, because the Poincar\'{e}
Lemma doesn't hold. One could imagine working with polar
coordinates $(r,\theta)$ in each fibre. Each component $h_{ij}$
of $h$ is periodic in $\theta$; but because $F$ is obtained by
integrating (twice), there is no guarantee that an $F$ can be
found which is periodic.

We can give conditions on the $h_{ij}$, expressed in terms of
$r$ and $\theta$, for the existence of a periodic $F$ as
follows. We have $y^1=r\cos\theta$, $y^2=r\sin\theta$, so the
$h_{ij}$ satisfy \begin{align*}
h_{11}\cos\theta+h_{12}\sin\theta&=0\\
h_{21}\cos\theta+h_{22}\sin\theta&=0,
\end{align*}
whence $h_{12}=-(h_{11}+h_{22})\sin\theta\cos\theta$. It will in
fact be most convenient to work in terms of the trace
$h_{11}+h_{22}$. Now the $h_{ij}$ are homogeneous of degree
$-1$, which means that $r(h_{11}+h_{22})$ is a function of
$\theta$ alone:\ we denote it by $\tau(\theta)$. It is of course
periodic. We have $h_{11}=\tau\sin^2\theta/r$,
$h_{22}=\tau\cos^2\theta/r$,
$h_{12}=-\tau\sin\theta\cos\theta/r$. Then
\begin{align*}
h_{11}dy^1+h_{12}dy^2&=-(\tau\sin\theta) d\theta\\
h_{21}dy^1+h_{22}dy^2&=(\tau\cos\theta) d\theta.
\end{align*}
Necessary and sufficient conditions for these two 1-forms to be exact are
\[
\int_0^{2\pi}(\tau\sin\theta) d\theta
=\int_0^{2\pi}(\tau\cos\theta) d\theta=0.
\]

This doesn't completely answer the question, however, because
there is another integration to carry out. It turns out to be
better to start from scratch. Let us compute the Hessian of a
positively-homogeneous function $F$ in polar coordinates. We may
set $F(r,\theta)=r\varphi(\theta)$, where $\varphi$ is periodic.
Then by straightforward calculations
\begin{align*}
h_{11}&=\frac{1}{r}(\varphi''+\varphi)\sin^2\theta\\
h_{22}&=\frac{1}{r}(\varphi''+\varphi)\cos^2\theta\\
h_{12}&=-\frac{1}{r}(\varphi''+\varphi)\sin\theta\cos\theta.
\end{align*}
So $\varphi''+\varphi=\tau$, and the question is whether this
equation has a periodic solution $\varphi$ for given periodic
$\tau$. We expect that $\int_0^{2\pi}(\tau\sin\theta)
d\theta=\int_0^{2\pi}(\tau\cos\theta) d\theta=0$ should be
necessary conditions for the existence of a periodic solution;
and indeed
\[
(\varphi''+\varphi)\sin\theta=\frac{d}{d\theta}(\varphi'\sin\theta-\varphi\cos\theta),\quad
(\varphi''+\varphi)\cos\theta=\frac{d}{d\theta}(\varphi'\cos\theta+\varphi\sin\theta)
\]
and the integrals of these functions evidently vanish if
$\varphi$ is periodic. The conditions
$\int_0^{2\pi}(\tau\sin\theta)
d\theta=\int_0^{2\pi}(\tau\cos\theta) d\theta=0$ are in fact
sufficient, as can be seen as follows. Define $u(\theta)$,
$v(\theta)$ by
\[
u(\theta)=-\int_0^\theta (\tau(\phi)\sin\phi)d\phi,\quad
v(\theta)=\int_0^\theta (\tau(\phi)\cos\phi)d\phi.
\]
Then $\varphi(\theta)=u(\theta)\cos\theta+v(\theta)\sin\theta$
is a particular solution of the equation
$\varphi''+\varphi=\tau$, as can easily be seen by direct
calculation (it is in fact the solution obtained by the method
of variation of parameters). Now
\[
u(\theta+2\pi)-u(\theta)=-\int_0^{2\pi} (\tau(\phi)\sin\phi)d\phi,\quad
v(\theta+2\pi)-v(\theta)=\int_0^{2\pi} (\tau(\phi)\cos\phi)d\phi;
\]
so if these integrals vanish then the equation has periodic
solutions. Notice that the addition to the particular solution
of a term $a\cos\theta+b\sin\theta$ (the complementary function)
corresponds merely to the addition of a term linear in $y$ to
$F$.

If $\tau=k_1\cos\theta+k_2\sin\theta$ there is no periodic $F$
whose Hessian is $h$; and indeed this is the key case, since for
any $\tau$ there are constants $k_1$ and $k_2$ (more exactly,
functions on $M$) such that $\tau-(k_1\cos\theta+k_2\sin\theta)$
does lead to a periodic $F$. We can nevertheless solve the
equation $\varphi''+\varphi=k_1\cos\theta+k_2\sin\theta$:\ a
particular solution is
\[
\varphi(\theta)=\onehalf(k_1\sin\theta-k_2\cos\theta)\theta.
\]
The expressions for the corresponding $h_{ij}$ are
\begin{align*}
h_{11}&=\frac{(k_1y^1+k_2y^2)(y^2)^2}{r^4}\\
h_{22}&=\frac{(k_1y^1+k_2y^2)(y^1)^2}{r^4}\\
h_{12}&=-\frac{(k_1y^1+k_2y^2)y^1y^2}{r^4}.
\end{align*}

When faced with a non-periodic $F$ one possible course of
action, which can certainly be carried out over a coordinate
neighbourhood $U\subset M$, is to replace the fibre by its
universal covering space (or in other words simply to ignore the
fact that $F$ is not periodic). It appears that in this way one
would obtain, on restricting to a level set of $F$, an example
of what Bryant calls a generalized Finsler structure (see
\cite{Bry}).

\section{Example}

The purpose of the example described below is to illustrate how
a spray may belong to the geodesic class of a globally-defined
pseudo-Finsler function, where the pseudo-Finsler function may
be made into a Finsler function only locally.

Consider the projective class of the spray
\[
\Gamma=u\vf{x}+v\vf{y}+w\vf{z}
+\sqrt{u^2+v^2+w^2}\left(-v\vf{u}+u\vf{v}\right)
\]
defined on $T^\circ\R^3$.  The geodesics of $\Gamma$ are spirals with
axis parallel to the $z$-axis, together with straight lines parallel
to the $z$-axis and circles in the planes $z=\mbox{constant}$.
To see this, note first that both $\sqrt{u^2+v^2}=\mu$ and $w$ are
constant; and therefore (or directly) $\sqrt{u^2+v^2+w^2}=\lambda$ is
also constant.  The geodesics are solutions of
\[
\ddot{x}=-\lambda\dot{y},\quad
\ddot{y}=\lambda\dot{x},\quad
\ddot{z}=0.
\]
Integrating the first two we get
\[
\dot{x}=-\lambda(y-\eta),\quad \dot{y}=\lambda(x-\xi)
\]
with constants $\xi$, $\eta$, whence
\[
(x-\xi)^2+(y-\eta)^2=(\mu/\lambda)^2.
\]
So the projections of the geodesics on the $xy$-plane are circles
of center $(\xi,\eta)$ and radius $r=\mu/\lambda$:\ note that
$0\leq r\leq 1$, the circle degenerating to a point when $r=0$.
The explicit parametrization of the geodesics is
\[
x(t)=\xi+r\cos(\lambda t+\vartheta),\quad
y(t)=\eta+r\sin(\lambda t+\vartheta),\quad
z(t)=wt+z_0,
\]
where $\xi,\eta,r,\lambda,\vartheta, w$ and $z_0$ are constants, with
$w^2=\lambda^2(1-r^2)$. So for $w/\lambda\neq0,\pm1$ the
geodesics are spirals, with axis the line parallel to the
$z$-axis through $(\xi,\eta,0)$. The case $r=0$ corresponds to
$w/\lambda=\pm1$ and the geodesics are straight lines parallel
to the $z$-axis (in both directions). The case $r=1$ ($w=0$)
gives circles of unit radius in the planes $z=z_0$.

In fact $\Gamma$ belongs to the geodesic class of the pseudo-Finsler function
\[
F(x,y,z,u,v,w)=\sqrt{u^2+v^2+w^2}+\onehalf yu-\onehalf xv.
\]
This is globally well defined but only locally a Finsler
function:\ it is positive only for $x^2+y^2<4$. It is globally
pseudo-Finsler, however. To obtain a Finsler function in a neighbourhood
of an arbitrary point $(x_0,y_0,z_0)$ we can make a simple
modification to
\[
\tilde{F}(x,y,z,u,v,w)=\sqrt{u^2+v^2+w^2}
+\onehalf(y-y_0)u-\onehalf(x-x_0)v;
\]
this is positive for $(x-x_0)^2+(y-y_0)^2<4$. Note that it
differs from $F$ by a total derivative.

The planes $z=\mbox{constant}$ have the property that a geodesic
initially tangent to such a plane (so that $w=0$ initially)
remains always in the plane:\ that is to say, such planes are
totally-geodesic submanifolds. The restriction of $\Gamma$ to
the submanifold $z=0$, $w=0$ of $T^\circ\R^3$ is the spray
\[
u\vf{x}+v\vf{y}-v\sqrt{u^2+v^2}\vf{u}+u\sqrt{u^2+v^2}\vf{v}
\]
of Shen's circle example from \cite{Shen}. We consider this as a
spray defined on $T^\circ\R^2$. It has for its geodesics all
circles in $\R^2$ of radius 1, traversed counter-clockwise.
Again, this spray is locally projectively metrizable. One local
Finsler function is the restriction of the one given for the
spiral example, namely
\[
F(x,y,u,v)=\sqrt{u^2+v^2}+\onehalf yu-\onehalf xv.
\]
Again, $F$ is only locally defined as a Finsler
function, though it is global as a pseudo-Finsler
function.

Finally we shall compare the geodesics of the Finsler function
\[
F(x,y)=\sqrt{g_{ij}(x)y^iy^j}+\beta_i(x)y^i=\alpha(x,y)+\beta(x,y)
\]
of Randers type (to which class the examples above belong) with
motion under the Lagrangian
\[
L(x,y)=\onehalf g_{ij}(x)y^iy^j+\beta_i(x)y^i.
\]
Here $g$ is a Riemannian metric. In two- or three-dimensional
flat space this is the Lagrangian for the motion of a classical
charged particle, of unit charge, in the magnetic field determined by
$d(\beta_idx^i)$.

We have
\[
\fpd{L}{y^i}=g_{ij}y^j+\beta_i=y_i+\beta_i.
\]
The energy $E_L$ is given by
\[
E_L=\Delta(L)-L=g_{ij}y^iy^j+\beta_iy^i-L=\onehalf
g_{ij}y^iy^j=\alpha.
\]
Thus motion under $L$ is with constant Riemannian speed. The
Euler-Lagrange field $\Gamma_L$ is determined by
\[
\Gamma_L(y_i)=\alpha\fpd{\alpha}{x^i}
-y^j\left(\fpd{\beta_i}{x^j}-\fpd{\beta_j}{x^i}\right).
\]
On the other hand
\[
\fpd{F}{y^i}=\frac{y_i}{\alpha}+\beta_i,
\]
and any geodesic spray $\Gamma_F$ of $F$ satisfies (but is
not determined by)
\[
\Gamma_F(y_i)=\Gamma_F(\alpha)\frac{y_i}{\alpha}+\alpha\fpd{\alpha}{x^i}
-\alpha y^j\left(\fpd{\beta_i}{x^j}-\fpd{\beta_j}{x^i}\right).
\]
We can fix the spray $\Gamma_F$ by the requirement that
$\Gamma_F(\alpha)=0$, so that motion is with constant Riemannian
speed (this is not of course the canonical spray of $F$). Then
for motion with unit Riemannian speed, motion under the
Lagrangian $L$ is the same as motion along the geodesics of the
fixed spray $\Gamma_F$. (For another speed one would have to
replace $\beta$ by a constant multiple of it.)

It has been claimed \cite{Tab} that the derivation of the equations of
motion from $F$ rather than $L$ is an example of Maupertuis's
principle in action. Be that as it may, this analysis does
reveal that the spiral example above can be related to the
motion of a charged particle in a constant magnetic field along
the $z$-axis. It is known that in this regime circular motion is
possible in planes perpendicular to the direction of the
magnetic field, the radius of such circles being a constant
whose value depends on the strength of the field and the charge
(amongst other things); it is known as the gyroradius, Larmor
radius or cyclotron radius. This is of course Shen's circle
example.

\section{Concluding remarks}

It is pointed out in \cite{Shen} that there cannot be a globally-defined
Finsler function for the circle example above, for the
following reason. In a Finsler space whose geodesic spray is
positively complete, that is, for which every geodesic is
defined on $[0,\infty)$, every pair of points in $M$ can be
joined by a geodesic. This is a conclusion of the Hopf-Rinow
Theorem of Finsler geometry. In the circle example the geodesics
are evidently positively complete, but equally evidently there
are pairs of points which cannot be joined by a geodesic. The
same is true of the spiral example.

This observation raises an interesting point about the different
levels at which global questions enter the problem. We leave
aside the two-dimensional case, which as we have seen is
atypical. The differential conditions on a multiplier stated in
Theorem \ref{mult} are merely the starting point for further
analysis which generally proceeds as follows. These differential
conditions are regarded as partial differential equations for
the unknowns $h_{ij}$. In principle these differential equations
generate integrability conditions, which are further conditions
on the coefficients $\Gamma^i$ and their derivatives. In
favourable cases, for example when the spray is isotropic, such
integrability conditions are satisfied, and one can assert the
local existence of a multiplier satisfying the differential
conditions. But there is no guarantee that such a multiplier is
even defined y-globally.

If there is in fact a y-global multiplier, and it is everywhere
positive quasi-definite, then there is always a Finsler function
which is local in $M$, that is, defined y-globally over an open
subset of $M$. There is a cohomological obstruction on $M$ to
combining such locally-defined Finsler functions; but even when
this obstruction vanishes the result is globally only a
pseudo-Finsler function in general, and can be turned into a genuine
Finsler function (by addition of a total derivative) only
locally. This is the situation exemplified in the previous section.

It seems that to make further progress one would have to impose
some conditions on the global properties of the sprays of the
projective class. The most important question is how one
incorporates the observation above about completeness and the
Hopf-Rinow property into the story --- a question, we suspect,
of some subtlety because ideally it should be answered in a
projectively-invariant manner.

\section*{Appendix A:\ some basic formulas}

Let $\Gamma$ be a spray
\[
\Gamma=y^i\vf{x^i}-2\Gamma^i\vf{y^i},
\]
with corresponding horizontal and vertical local vector field basis
$\{H_i,V_i\}$ where
\[
H_i=\vf{x^i}-\Gamma^j_i\vf{y^j},\quad
\Gamma^i_j=\fpd{\Gamma^i}{y^j}.
\]
We denote by $R^j_i$ the components of the Riemann curvature or
Jacobi endomorphism of $\Gamma$. We have the following well-known formulas, which serve to define $R^j_i$:
\[
[\Gamma,H_i]=\Gamma^j_iH_j+R^j_iV_j,\quad
[\Gamma,V_i]=-H_i+\Gamma_i^jV_j,
\]
Of course $[V_i,V_j]=0$; we have
\[
[V_i,H_j]=-\conn kijV_k=[V_j,H_i],\quad \conn kij
=\fpd{\Gamma^k_i}{y^j}.
\]
We set
\[
R^i_{jk}=\frac13\left(\fpd{R^i_j}{y^k}-\fpd{R^i_k}{y^j}\right);
\]
then $R^i_j=R^i_{jk}y^k$, and the bracket of horizontal fields is
determined by
\[
[H_j,H_k]=-R^i_{jk}\fpd{}{y^i}.
\]

Associated with any spray there is an operator $\nabla$ on
tensors along $\tau$ which is a form of covariant derivative,
and indeed is often called the dynamical covariant derivative.
We need it mainly for its action on tensors $h$ with
components $h_{ij}$, when
\[
(\nabla h)_{ij}=\Gamma(h_{ij})-\Gamma_i^kh_{kj}-\Gamma_j^kh_{ik}.
\]

If $\Gamma$ is a spray,  another spray $\tilde{\Gamma}$ is
projectively equivalent to $\Gamma$ if there is a function $P$ on
$\TMO$, necessarily of homogeneity degree 1, such that
$\tilde{\Gamma}=\Gamma-2P\Delta$. Then
\[
\tilde{\Gamma}^i_j=\Gamma^i_j+P\delta^i_j+P_jy^i,
\quad P_j=\fpd{P}{y^j}.
\]

\section*{Appendix B:\ geodesic convexity and good open coverings}

In Section \ref{glob} we stated that every manifold admits good open
coverings by coordinate patches. The usual argument for this
uses Whitehead's results in \cite{CR} on the existence of
geodesically convex sets, together with the fact that on any
paracompact manifold one can construct a Riemannian metric,
which of course provides a source of geodesics.

Though it is convenient, it is not actually necessary to appeal
to the existence of a Riemannian metric. In the first place the
results on convexity in \cite{CR} apply to the geodesics of any
affine connection. What seems not to be so well known is that
they also apply to the geodesics of any spray, when due account
is taken of the fact that in this case the geodesics will be not
in general be reversible:\ that is to say, when one remembers
always to speak of the geodesic path from $x_1$ to $x_2$ (which
may not be the same as the geodesic path from $x_2$ to $x_1$).
Whitehead himself acknowledged that the results of \cite{CR} can
be extended in this way, in an addendum to that paper,
\cite{CRadd}, published the following year. Since in our work we
always have a spray at our disposal it seems natural to use this
second version of the convexity result. The fact that
geodesically convex regions exist for spray spaces is of
interest more widely in spray and Finsler geometry than just in
relation to the construction of good open coverings. Whitehead's
addendum \cite{CRadd} appears to be far less well known than the
original paper \cite{CR}; and it seems fair to say that in it
Whitehead was not as careful about the matter of the
non-reversiblity of geodesics of sprays as he might have been.
For these reasons we have thought it worthwhile to outline how
the result is proved.

Let us be precise about what is to be proved. Following
Whitehead we shall use the term {\it open region\/} for an open
subset of a coordinate patch of a manifold $M$, {\it closed
region\/} for the closure of an open region, and {\it region\/}
for either. A region $C$ is (geodesically) {\it convex\/} with
respect to a spray $\Gamma$ on $\TMO$ if for every $x_1,x_2\in
C$ there is at least one geodesic path of $\Gamma$ from $x_1$ to
$x_2$ lying entirely within $C$. A region $C$ is (geodesically)
{\it simple\/} if there is at most one geodesic path of $\Gamma$
from $x_1$ to $x_2$ lying entirely within $C$. Whitehead proves
(for affine sprays in \cite{CR}, and for sprays in general in
\cite{CRadd}) that for any manifold $M$ equipped with a spray
$\Gamma$, ``any point in $M$ is contained in a simple, convex
region which can be made as small as we please''.

The basic analytical tool is Picard's Theorem on the existence
and uniqueness of solutions of two-point or boundary-value
problems for systems of second-order ordinary differential
equations. Using this theorem Whitehead is able to show that for
any $x\in M$ there is a region $C_1$ containing $x$ such that
for $x_1,x_2\in C_1$ there is a geodesic
$s\mapsto\gamma(x_1,x_2,s)$ with $\gamma(x_1,x_2,0)=x_1$,
$\gamma(x_1,x_2,1)=x_2$, and $\gamma(x_1,x_2,s)\in C_1$ for
$0\leq s\leq1$; moreover $\gamma$ is continuous in all of its
arguments. Thus $C_1$ is convex. Whitehead shows further that
there is a second region $C_2$ with $x\in C_2\subset C_1$ which
is simple. These results apply to any spray, not
just to affine ones; it was Whitehead's realization of this
point that led him to write his addendum.

Now any subregion of a simple region is simple. This is not true
of a convex region however:\ a subregion of a convex region need
not be convex. But suppose that we can find a subregion of a
convex region $C_1$ whose interior is an open connected set $C$
containing $x$, whose closure is compact, and whose boundary $B$
is a smooth hypersurface (codimension 1 submanifold) with the
property that any geodesic tangent to $B$ passes outside
$\bar{C}=C\cup B$ at least locally:\ that is to say, that if
$\gamma$ is a geodesic with $\gamma(0)\in B$ and
$\dot{\gamma}(0)$ tangent to $B$, then $\gamma(s)\notin \bar{C}$
for $s$ in some open interval $(0,t)$, $t>0$. In this case a
geodesic cannot touch $B$ while remaining otherwise in $C$.
Following a line of argument from \cite{CR}, using this
observation, we show that $C$ must be convex.

Consider the set of points $(a,b)\in C\times C$ such that there
is a geodesic in $C$ from $a$ to $b$:\ call it $G$. Firstly, $G$
is not empty:\ for any $a\in C$ and for any geodesic $\gamma$
starting at $a$, $\gamma(s)$ must lie in $C$ for all $s$ in some
interval $[0,t)$ for $t$ sufficiently small. Secondly, $G$ is open by
continuity of $\gamma(a,b,s)$. Thirdly, $G$ is closed (as a
subset of $C\times C$). Consider a point $(a,b)$ of $C\times C$
which lies in the relative closure of $G$; that is, $(a,b)$ is
the limit of a sequence $(a_n,b_n)$ of points of $G$. Each
geodesic $s\mapsto\gamma(a_n,b_n,s)$ lies in $C$ for $0\leq
s\leq 1$. There is a geodesic $s\mapsto\gamma(a,b,s)$ from $a$
to $b$, but we know only that it lies in $C_1$ (though of course
its initial and final points $a$ and $b$ are in $C$). Now for
each $s$ with $0\leq s\leq 1$, $\gamma(a,b,s)$ is the limit of
the sequence $\gamma(a_n,b_n,s)$, so $\gamma(a,b,s)$ certainly
lies in $\bar{C}$ for all $s$. But the geodesic
$s\mapsto\gamma(a,b,s)$ lies in $C$ initially, and can neither
meet $B$ transversely, nor be tangent to it, without
subsequently passing out of $\bar{C}$. So $\gamma(a,b,s)$ lies
in $C$ for all $s$, $0\leq s\leq1$, and $(a,b)\in G$. Since $C$
is connected by assumption, so is $C\times C$; $G$ is a
non-empty subset of $C\times C$ which is both open and closed,
and so $G=C\times C$.

So if we can find a region $C$ with these properties in $C_2$ it
will be convex and simple. In fact the open Euclidean coordinate
ball $\{(x^i):\delta_{ij}x^ix^j<r^2\}$ for any sufficiently
small positive $r$ will do for $C$ (we take $x$ as origin of
coordinates). Here $B$ is the Euclidean coordinate sphere of
radius $r$ of course. For any positive $r$ consider the function
$V_r(x^i)=\delta_{ij}x^ix^j-r^2$. Consider a geodesic $\gamma$
such that $V_r(\gamma(0))=0$ and
\[
\frac{d}{ds}\Big(V_r\circ\gamma\Big)_{s=0}
=2\delta_{ij}\gamma^i(0)\dot{\gamma}^j(0)=0;
\]
it is tangent to $B$ at $\gamma(0)\in B$. We have
\[
\frac{d^2}{ds^2}\Big(V_r\circ\gamma\Big)_{s=0}
=2\delta_{ij}\dot{\gamma}^i(0)\dot{\gamma}^j(0))
-4\delta_{ij}\gamma^i(0)\Gamma^j(\gamma(0),\dot{\gamma}(0)).
\]
We shall show that for $r$ sufficiently small this is positive,
which means that $V_r(\gamma(s))>0$ for $s$ in some open
interval about 0, but $s\neq0$, from which it will follow that
any geodesic tangent to $B$ locally lies outside $\bar{C}$.
Without essential loss of generality we may assume that the
first term is 2. For any $r_0>0$ the set
$\{(x,y):\delta_{ij}x^ix^j\leq r_0^2,\,\delta_{ij}y^iy^j=1\}$ is
compact, so there is $K>0$ such that $|\Gamma^i(x,y)|<K$ for all
$i$ and all $(x,y)$ in that set. Take $r\leq r_0$. Then
$|\gamma^i(0)|\leq r$ and
$|\Gamma^j(\gamma(0),\dot{\gamma}(0))|<K$, whence
\[
|\delta_{ij}\gamma^i(0)\Gamma^j(\gamma(0),\dot{\gamma}(0))|<nrK.
\]
So provided that $r<1/(2nK)$,
\[
\frac{d^2}{ds^2}\Big(V_r\circ\gamma\Big)_{s=0}>0
\]
as required. This argument is essentially the same as the one
given by Whitehead in \cite{CRadd}. It actually shows that the
final inequality holds for all geodesics, not just those tangent
to $B$; and in fact Whitehead gives this as his requirement for
$C$ to be convex, though this seems to us to be a stronger
condition than is really necessary. Be that as it may, we have
shown that for every sufficiently small $r$ the Euclidean
coordinate ball of radius $r$ about $x$ is a simple, convex
region.

A convex region $C$ is path connected and so connected. It is
also contractible on any point $x_0\in C$, for the map
$h:C\times [0,1]:h(x,s)=\gamma(x,x_0,s)$ is a continuous
homotopy of the identity with the constant map $x\mapsto x_0$.
Moreover the intersection of simple, convex sets is simple,
convex. Thus a covering by open simple, convex regions is a good
open covering.

\subsubsection*{Acknowledgements}
The first author is a Guest Professor at Ghent University:\ he
is grateful to the Department of Mathematics for its
hospitality. The second author is a Postdoctoral Fellow of the
Research Foundation -- Flanders (FWO). The third author
acknowledges the support of grant no.\ 201/09/0981 for Global
Analysis and its Applications from the Czech Science Foundation.

This work is part of the {\sc irses} project {\sc geomech} (nr.\
246981) within the 7th European Community Framework Programme.

\subsubsection*{Address for correspondence}
\noindent
M.\ Crampin, 65 Mount Pleasant, Aspley Guise, Beds MK17~8JX, UK\\
m.crampin@btinternet.com


\begin{thebibliography}{99}

\bibitem{AP}
J.\,C.\ \'{A}lvarez Paiva, Symplectic geometry and Hilbert's
fourth problem {\it J.\ Diff.\ Geom.} {\bf 69} (2005) 353--378.

\bibitem{APB}
J.\,C.\ \'{A}lvarez Paiva and G.\ Berck, Finsler surfaces with
prescribed geodesics, arXiv:1002.0243.

\bibitem{Bacso}
S.\ B\'{a}cs\'{o} and Z.\ Szilasi, On the projective theory of
sprays {\it Acta Math.\ Acad.\ Paed.\
Ny\'{\i}regyh\'{a}ziensis\/} {\bf 26} (2010) 171--207.

\bibitem{Bao}
D.\ Bao, S.-S.\ Chern and Z.\ Shen, {\it An Introduction to
Riemann-Finsler Geometry\/} Springer (2000).

\bibitem{Bry}
R.\,L.\ Bryant, Projectively flat Finsler 2-spheres of constant curvature
{\it Selecta Math.\ (N.\,S.)\/} {\bf 3} (1997) 161--203.

\bibitem{BM}
I.\ Bucataru and Z.\ Muzsnay, Projective metrizability and
formal integrability {\it SIGMA\/} {\bf 7} (2011) 114.

\bibitem{hil4}
M.\ Crampin, Some remarks on the Finslerian version of Hilbert's
fourth problem {\it Houston J.\ Math.} {\bf 37} (2011) 369--391.

\bibitem{Japan}
M.\ Crampin and D.\,J.\ Saunders, Path geometries and almost
Grassmann structures {\it Adv.\ Stud.\ Pure Math.}\ {\bf 48} (2007)
225--261.

\bibitem{Doug}
J.\ Douglas, Solution of the inverse problem of the calculus of
variations {\it Trans.\ Amer.\ Math.\ Soc.} {\bf 50} (1941)
71--128.

\bibitem{Hen}
M.\ Henneaux, On the inverse problem of the calculus of
variations {\it J.\ Phys.\ A:\ Math.\ Gen.} {\bf 15} (1982)
L93--96.

\bibitem{KrupkaS}
D.\ Krupka and A.\,E.\ Sattarov, The inverse problem of the
calculus of variations for Finsler structures {\it Math.\
Slovaca\/} {\bf 35} (1985) 217--222.

\bibitem{KP}
O.\ Krupkov\'{a} and G.\,E.\ Prince, Second order ordinary
differential equations in jet bundles and the inverse problem of
the calculus of variations, in {\it Handbook of Global
Analysis\/} (D.\ Krupka and D.\,J.\ Saunders, ed.) Elsevier
(2008) 837--904.

\bibitem{Lov}
R.\,L.\ Lovas, A note on Finsler-Minkowski norms {\it Houston
J.\ Math.} {\bf 33} (2007) 701--707.

\bibitem{Sar}
W.\ Sarlet, The Helmholtz conditions revisited. A new approach
to the inverse problem of Lagrangian dynamics {\it J.\ Phys.\
A:\ Math.\ Gen.} {\bf 15} (1982) 1503--1517.

\bibitem{Sar2}
W.\ Sarlet, Linear connections along the tangent bundle
projection, in {\it Variations, Geometry and Physics\/} (O.\
Krupkov\'{a} and D.\,J.\ Saunders, ed.) Nova Science (2009)
315--340.

\bibitem{Shen}
Z.\ Shen, {\it Differential Geometry of Spray and Finsler Spaces\/} Kluwer
(2001).

\bibitem{Szi}
J.\ Szilasi, Calculus along the tangent bundle projection and
projective metrizability, in {\it Proceedings of the 10th
International Conference on Differential Geometry and its
Applications, Olomouc, Czech Republic, 2007\/} (O.\ Kowalski,
D.\ Krupka, O.\ Krupkov\'{a} and J.\ Slov\'{a}k, ed.) World
Scientific (2008) 527--546.

\bibitem{SzVa}
J.\ Szilasi and Sz.\ Vattam\'{a}ny, On the Finsler metrizabilities of spray
manifolds, {\it Period.\ Math.\ Hungarica\/} {\bf 44} (2002) 81--100.

\bibitem{Tab}
S.\ Tabachnikov, Remarks on magnetic flows and magnetic
billiards, Finsler metrics, and a magnetic analogue of Hilbert's
fourth problem, in {\it Modern Dynamical Systems and
Applications\/} (M.\ Brin, B.\ Hasselblatt, and Y.\ Pesin, ed.),
Cambridge University Press (2004) 23--252.

\bibitem{Warner}
F.\,W.\ Warner, {\it Foundations of Differentiable Manifolds and Lie Groups\/}
Scott, Foresmann (1971).

\bibitem{CR}
J.\,H.\,C.\ Whitehead, Convex regions in the geometry of paths
{\it Quart.\ J.\ Math.} {\bf 3} (1932) 33--42.

\bibitem{CRadd}
J.\,H.\,C.\ Whitehead, Convex regions in the geometry of paths -- addendum
{\it Quart.\ J.\ Math.} {\bf 4} (1933) 226--227.



\end{thebibliography}
\end{document}